# Introducing an old calculating instrument in a new technologies environment: a praxeological analysis of students' tasks using different registers

**Caroline POISARD**

*Laboratoire du CREAD*
*UBO, ESPE de Bretagne, site de Quimper*
*France*
*caroline.poisard@espe-bretagne.fr*

**ABSTRACT**
*The Chinese abacus is the resource presented in this paper, to teach and learn number sense and place-value system at primary level. The Chinese abacus can be material, virtual (software) or drawn on a worksheet. We present three tasks and analyse them in term of techniques and relative knowledge. We show how these tasks can be solved by students in different registers (material, software, paper-and-pencil, fingers, oral) which is important for both students' understanding and teachers' activity.*

**RÉSUMÉ**
*Le boulier chinois est la ressource présentée dans cet article, pour enseigner et apprendre la construction du nombre et le système de numération décimal à l'école. Le boulier chinois peut être matériel, virtuel (logiciel) ou dessiné sur une feuille. Nous présentons trois tâches et les analysons en terme de techniques et connaissances sous-jacentes. Nous montrons comment ces tâches peuvent être résolues par des élèves dans différents registres (matériel, logiciel, papier-crayon, mains, oral) ce qui est important pour la compréhension des élèves et également pour l'activité des professeurs.*

**MOTS-CLÉS**
*Ressources matérielles et virtuelles, praxéologie, technique, technologie, registre, construction du nombre, numération décimale, boulier chinois*

**INTRODUCTION**

Nowadays, in most classrooms in France, not only material resources are available for teachers and students, but as well virtual ones. Computers are used to prepare their class by teachers, students can have access to computers or tablets in the classroom or in a computer room, teachers in the classroom can use a video-projector, an interactive whiteboard (IWB) with different softwares, etc. The introduction of symbolic calculators at secondary level in the 90's gave rise to some very interesting works in mathematics education (Guin, Ruthven & Trouche 2005). It seems to us these analysis are also appropriate to look at some more *traditional* resources. In our study, teachers use the Chinese abacus (suan-pan) for teaching and learning number sense, place-value system and calculation at primary level in France. The fact that we precise it is in France is important: the abacus is not introduced as it can be in Asia, to become expert and be able to calculate easily and fast. The abacus in studied to introduce a new register to work on number sense and place-value system. We show here that it is as well a good opportunity for inquiry based learning. A software is also used -that we call the virtual abacus- and articulated to the material abacus.

Numbers and arithmetic is a rich mathematical area for education. A recent ICMI study (ICMI Study 23, Sun, Berinderjeet and Novotna 2015) presents studies on the teaching and the learning of whole numbers at primary level. Number sense is a difficult notion to define. For example, Baccaglini-Frank and Maracci (2015) give a good overview of the question in the fields of mathematics education and cognitive psychology as well. In this paper, by *number sense* we mean: the distinction between the ordinality and the cardinality of numbers[1], the difference between value and quantity, and the decomposition of numbers.

The first part of our paper presents the theoretical frameworks we retain for this study and our methodological choices. The second part identify the important particularities of the Chinese abacus, including the particularities of the software. With some examples of teachers' discourse about their experience in class we show how register is an underlying notion. The third part is the praxeological analysis of three tasks given to students in different registers. The initial questions raised are: Considering some tasks using the Chinese abacus, how can we describe the knowledge used by students? How the register can (or can not) influence students mathematical activity? Our research questions are exposed at the end of the first part.

## 1. THEORETICAL REFERENCES AND METHODOLOGICAL CHOICES

*Resources*
The term of resource is seen in a large meaning, referring to Adler's work:
"The common-sense notion of resources in and for education is resource as a material object, and lack of resources usually refers to shortages of textbooks and other learning materials. It is possible to think about resource as the verb *re-source*, to source again or differently." (Adler, 2000, p.207).

So, the main point for a resource is the possibility to resource the teaching and the learning of mathematics. For Adler, a resource can be -at least- material (paper, blackboard, textbook, computer, etc.), human (a discussion with a colleague, students' work, etc.), or cultural (different languages used in the class by teachers and/or students, etc.). Part of our work is about *documentational resources* (Gueudet and Trouche, 2009) but this is not the focus we have chosen to develop here.

We introduced the distinction between *material* and *virtual resource* (Poisard, Gueudet and Bueno-Ravel, 2011). It seems to us that, in education, new technologies are not replacing some other old artifacts. We can see a reorganisation of the articulation of different resources for teachers using new technologies. At primary level, by analysing the appropriation of the Chinese abacus by teachers, we have shown that to teach number sense and place-value system, teachers articulate the material Chinese abacus, a virtual one (software used with an interactive whiteboard IWB), and paper-and-pencil activities on worksheets. Maschietto and Trouche (2010) discussed as well this idea when old and new technologies are present in the same class situation, they give an analysis in term of instrumentation and orchestration. Our work is connected with other research, using some different frameworks, that introduced some terms that can have similarities with the notion of resources. For examples, in their paper on "manipulatives in mathematics education", Bartolini and Martignone (2014) refer to *concrete* and *virtual manipulatives* but also to *historico-cultural* and *artificial manipulatives* (designed by educators). In their work referring to the instrumental approach, Maschietto and Soury-Lavergne (2013) make a distinction between *material* and *digital artifacts*.

---

1  The ordinality is the position of a number into the natural number sequence (ex: the third bead) and the cardinality is a quantity (ex: three beads).

*Registers*

Part of our work here is to analyse some tasks given to students and to show how the different resources can be articulated by teachers. We look at teachers' work at a small scale, and we refer to the notion of *register* introduced by Duval (1996, 2006) in his cognitive approach. Here registers are resources, they are included in a larger range of resources available for a teacher. Of course, all resources are not a representation register. Duval's work is centred on the understanding of mathematical difficulties of students. A classification of *registers of semiotic representation* is proposed by Duval to analyse mathematical activity. First what is a *representation*? For Duval:

"Representations can be individuals' beliefs, conceptions or misconceptions to which one gets access through the individuals' verbal or schematic productions. […] But representations can also be signs and their complex associations, which are produced according to rules and which allow the description of a system, a process, a set of phenomena. There the semiotic representations, including any language, appear as common tools for producing new knowledge and not only for communicating any particular mental representation." (Duval, 2006, p.104).

What is important for the mathematical activity is the *semiotic representations*. For Duval, *semiotic representations* is central to define the mathematical process:

"No kind of mathematical processing can be performed without using a semiotic system of representation, because mathematical processing always involves *substituting some semiotic representation for another*. The part that signs play in mathematics is not to be substituted for objects but for other signs! What matters is not representations but their transformation." (Duval, 2006, p.107).

And with the notion of transformation, Duval defines the term *register*. A *representation register* is a semiotic system that allows *transformation* of representations, so not all semiotic representations are registers. For example, a verbal form, a numerical expression, a formal notation, a figure are some representation registers. More precisely, Duval considers two very different types of transformations: the *treatment* (that stays in the same register) and the *conversion* (where different registers are involved). For example, to solve an equation, one can stay in the register of notations, using symbols and mathematical rule to solve it, this will be a treatment transformation. But a conversion transformation coordinates at least two registers. For example, the task: "from a mathematical expression given in words, write it with symbols using =, >, <" means to coordinate a words expression with a symbolic one. A question important for us: can a material object be a semiotic representation? In this paper, Duval mentions as a note that:

"It is only from a strict formal point of view that semiotic representations can be taken as concrete objects (Duval, 1998, pp. 160–163)." (Duval, 2006, p.129).

We show in this paper that the material and virtual abacus can be considered as semiotic representation and articulated between different registers to explain the mathematical activity. For Duval, this is the important part of the mathematical activity of students, the hypothesis given is that "comprehension in mathematics assumes the coordination of at least two registers of semiotic representation" (Duval, 2006, p.115). For both theoretical frameworks, the cognitive and the documentational approach, an important feature to look at is the coordination, the articulation of registers, of resources. We show in the part 3 the articulation of different registers: verbal (oral and written languages), material abacus, virtual abacus, fingers and paper-and-pencil worksheets.

*Tasks and techniques*

The *anthropological theory of the didactic* (ATD, Chevallard 1992a, 1992b, 2006) develops several main concepts as the *didactic transposition* in a given *institution*, the *dialectic of medias and milieu,* and the *praxeological analysis*. We focus here on the praxeological analysis introduced to analyse the mathematical knowledge taught in institutions (Chevallard, 1999, 2007):

"Essentially praxeology is made of two parts, the *praxis* part and the *logos* part. […] The *praxis* part is the union of a *type of tasks* (such as solving quadratic equations, blowing one's nose, composing a fugue, for example) and a *technique* – way of doing – which purportedly allows one to carry out at least *some* tasks of the given type - those in the "scope" of the technique. The *logos* part is the union of a whole set of notions and arguments arranged into a more or less rational "discourse" (*logos*), the so called *technology* of the technique, which is intended to provide justification for the technique – why does it work (at least sometimes), where does its effectiveness come from?, etc. - and a more abstract set of concepts and arguments arranged into a more general "discourse", the praxeology's *theory*, supposed to justify the technology itself." (Chevallard 2007, p. 131).

The praxeological organisation has a practical part called the technique (*praxis*) and a discourse on the technique called the technology (*logos*). So, to analyse a given *type of tasks*, some *techniques* (at least one) are presented, and for each technique the associated *technology* (the mathematical knowledge) is developed. Finally, this is linked to a more general mathematical *theory*. For example, in geometry, it is central to mention if working in the euclidean geometry or not, to analyse some type of tasks.

*Related research using these theoretical frameworks*
In a previous research, we have used the articulation of the cognitive and the anthropological approach to describe some tasks depending on the register (Poisard 2005a, 2005b). We had chosen three tasks: to set a digit, to set a number and to add two numbers (with a carried-number), and two registers: "paper-and-pencil and oral" and "material, gesture with the Chinese abacus". Some following studies reinvested this articulation, including the virtual abacus (Poisard, Gueudet & Bueno-Ravel, 2011, Riou-Azou, 2014, Bueno-Ravel and Harent, 2016).

The introduction in classrooms of new technologies has started to question researchers in mathematics education since many years. In our point of view, these investigations are very interesting and can help understand the teaching and the learning using any resource, material or virtual. The book coordinated by Guin, Ruthven and Trouche (2005), shows how the introduction of symbolic calculators in the classroom leads to some "didactical challenges". In particular, referring to the ATD, Lagrange (2005) analyses some tasks given to students in terms of techniques and presents "a variety of new techniques […] related to paper-and-pencil techniques" (p.113). He introduces the notion of *instrumented technique* defined in the following chapter by Trouche (2005). The notion of instrument is related to the instrumental approach, an instrument is an artifact used by someone with some action schemes (instrument=artifact+schemes, Verillon and Rabardel, 1995). In this environment with calculators, the particular schemes are called *instrumented action schemes*. The symbolic calculator is not only a "technical" tool but an instrument for the mathematical activity. Analysing students' activity, Trouche points "the difficulty of moving from one register to another one [...]." (Trouche, 2005, p143). In their work, the instrument studied is the symbolic calculator but could be any instrument.

In an overview of the ATD, Artigue (2009) analyses the relationships and connections between theoretical frameworks in the case of the ATD. In this review, the articulation of the ATD with the instrumental approach (quoting in particular Lagrange and Trouche's work) and the semiotic dimension (quoting Duval in particular) are mentioned. The articulation of praxeology and representation registers is also developed in Block, Nikolantonakis, and Vivier (2012): the importance of representation registers is analysed with different type of tasks, and looking at three different countries combined influences of the register are shown.

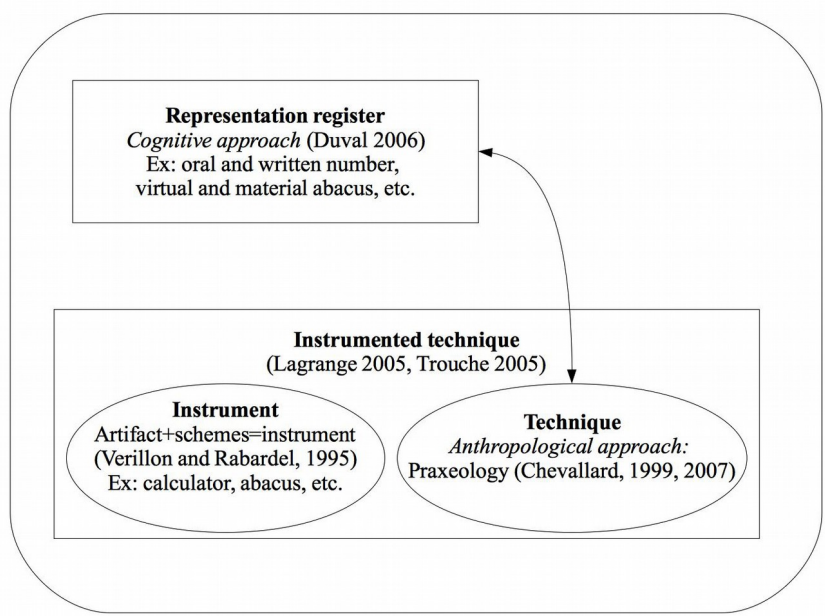

**FIGURE 1:** Articulation of different frameworks[2]

*Methodological choices*
This study is included in a larger project supported by the French Ministry of Education to develop resources for the teaching and learning of mathematics at primary level, including pre-school one. This national project has now come to an end, with the delivery of two training paths that are available for primary teachers training in mathematics, for both in-service and pre-service education (this national platform is called M@gistère). During those years, we worked within a research group composed by primary teachers and researchers for a collaborative work. For the Chinese abacus, at first a presentation was made of both the virtual and the material one, and teachers retain it to integrate in their progression to teach number sense and place-value system. We followed closely at least five primary teachers during five years (2010/15) from year 1 (5-6 years old) to year 6 (10-11 years old). That means they were regularly participating to the working group meetings to explain how the integration of the abacus was made in their class. The discussions were about: their choices (lesson plans), students' activity, the integrated resources (computers, pencil-and-paper activity, etc.), the technology available and the material constraints, the produced resources, etc. For the new resources produced by teachers, they were discussed in the group and ameliorated to produce the ones now available in the training path "The Chinese abacus at school" delivered in June 2015 (for a presentation see Poisard and al. 2016).

The data we analyse are:
- Notes from the working group meetings
- Observations and recordings of lessons
- Students' work (paper-and-pencil worksheets, recordings)
- Teachers' work (lesson plans, created resources)
- Interviews with teachers and students (5 to 11 years old)

To analyse these data, some recordings of interviews and class observations were transcribed. This paper is a synthesis of the research and all these data are combined for our analysis.

Our analysis answers the following questions: What are the registers we can identify? What is the mathematical activity depending the register? For a given task, can we identify different

---

2  For this figure, there is no scale between the frameworks. The aim is to show interactions between these frameworks.

techniques? Is it dependent of the register? Why?

## 2. THE CHINESE ABACUS: IDENTIFYING REGISTERS

*The virtual and material Chinese abacus*
Like the material one, the virtual abacus has two parts: a lower part and an upper part. Each rod corresponds to a rank of the place-value system: units, tens, hundreds, etc. (from the right to the left). On the lower part, the beads represent a value of one (ten, one hundred, etc. depending on the rod) and are called *one-unit counters*. On the upper part, they represent a value of five (fifty, five hundreds, etc. depending on the rod) and are called *five-unit counters*. The software we use (from Sésamath) is available online[3]. The software has three icons:
- "see number" that can display (or not) the number written in digits
- "set to zero"
- "positioning" that allows from any inscription to set the *economical inscription* (Poisard, 2005b) that means moving the less number of beads possible. For example, to set 25 we have several possibilities: to activate 2 one-unit counters in the tens and 1 five-unit counter in the units which is the *economical inscription* (A); but it is also possible to activate 2 one-unit counters in the tens and 5 one-unit counters in the units (B) or also; to activate 1 one-unit counter in the tens and in the units: 5 one-unit counters and 2 five-unit counters (C) (Figure 2). From the two last propositions, using the icon "positioning" will display the economical inscription of 25.
    So, on the figure 2, we have three different decompositions of 25:
- Inscription A: economical inscription: 25=20+5=(10+10)+5
- Inscription B: 25=20+5=(10+10)+(1+1+1+1+1)
- Inscription C: 25=10+15=(10)+(5+5+1+1+1+1+1).
    These three inscriptions are three different techniques to set 25 on the Chinese abacus. For each technique, there is a specific mathematical knowledge associated concerning the decomposition of 25. We develop the different possibilities to set a number and the mathematical knowledge associated in part 3 with the analysis of different tasks depending on the register.

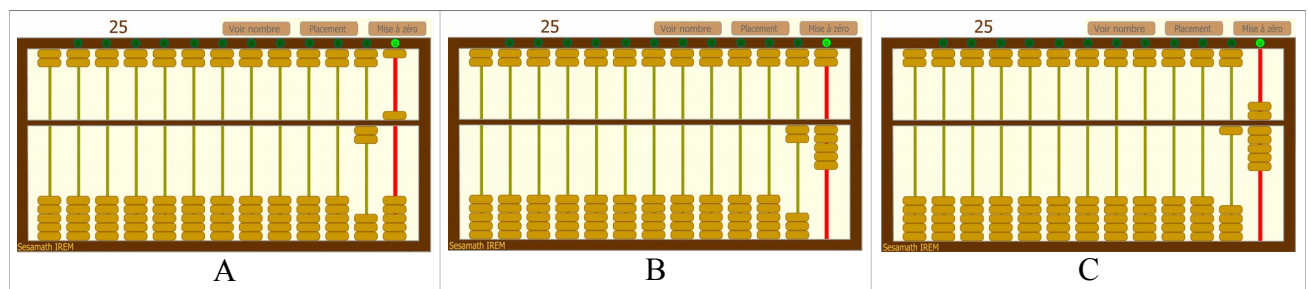

**FIGURE 2:** Three inscriptions of 25: economical inscription (left) and the two other ones

    Rau et al (2016) propose a model of the cognitive process of Chinese abacus arithmetic, describing three methods for solving arithmetic method problems: retrieval method, procedure method and mental arithmetic method. In this paper, we focus on number sense and place-value system (not on arithmetic problems). The two main tasks we identify to work on number sense and place-value system with the Chinese abacus are: "*to set a number*" and "*to read a number*". The feedback given by the software is the possibility to display a number written in digits. That means a student can use the icon "see number" to verify her/his answer and correct it if needed. But there is no recording of students' activity. That means that most teachers ask to answer on a paper-and-pencil worksheets. Teachers can also ask a student to perform the task using the abacus with an

---
3   The Sésamath Chinese abacus online: http://cii.sesamath.net/lille/exos_boulier/boulier.swf . The pictures in this paper are taken from this software.

IWB and discuss it.

*Teachers' discourse and the notion of register*
In teachers' discourse, we observe that there is a mention of different registers used with students. We refer here to two interviews with primary teachers: the first one was made with Ester (CM2, year 6) and Deborah (CP, year 2) (Poisard and all 2016) and the second interview made with Rose (CM1, year 5)[4]. Ester and Deborah had used the Chinese abacus for several years, and it is the first year for Rose. The three of them have more than fifteen years of experience as primary teachers (including preschool[5]) and teach as well trainers at University in-service and pre-service primary teacher education.

We transcribed these interviews and give here some extracts (translated from French). The interviews have three main directions concerning:
- the mathematical knowledge: About the Chinese abacus, what contribution did you notice in the class? What is the interest to use the abacus in class?
- students' activity: What learning seems to have been made? What are the procedures implemented? What are the difficulties and mistakes met with?
- teacher's work: What is teacher's role for these sessions? How to manage the class activity? How to support students' activity?

In their interviews, Ester, Deborah and Rose point: the differentiation, the students' autonomy, the motivation, the verbalisation of procedures, and their vision of mathematics and their classroom practice. In term of registers, we identify in teachers' spontaneous discourse: the material and the virtual abacus, the paper-an-pencil worksheet, and the verbalisation of numbers. Moreover, we have observed in their classrooms the use of fingers to set numbers and calculate.

Ester explains that the paper-and-pencil worksheets are often used to complete the virtual abacus that allows a self-correction for students: "We designed, little by little, some worksheets to make students more and more autonomous, so they can improve, work at their own rhythm [...]". Ester carries on with the interest of using the abacus to differentiate students' work: "Thanks to this tool, we can easily allow students to acquire number sense and calculation, with worksheets of differentiation so that students can effectively improve." For Rose, the variety of formats for students is very important, the material and virtual abaci are some more formats integrated in the classroom: "Here, students really manipulate numbers, when counting they move the beads or point them. There is a kinaesthetic side that is involved here, particularly for our students that need to move, to touch."

For Ester, the second interest is about computer skills: "With the virtual abacus, we associate computer skills with mathematical skills, and this is not neutral for students' motivation!". Rose notices as well that students like the abacus sessions: they are motivated and explain that with two arguments. Students like to work in groups and on computers. "Most of the sessions were in group, homogeneous groups with a resource person per group. For a tutelage work, with discussions and the virtual abacus to verify answers. [...] If students have to choose between the material abacus and the virtual one to manipulate; undeniably they choose the virtual one!"

Moreover, Deborah underlines that for students, sessions with the abacus allow to work on verbalisation, in particular on procedures explanation: "What is really obvious is students' verbalisation, they have to use the appropriate vocabulary: one-unit counter, five-unit counter. I ask them to make explicit their procedures […] after, they are used to be in this approach, to really make explicit their procedure, that can be an expert one or not. And it is what is interesting […]. Everyone has the benefit from this discussion!". Indeed, the integration of the Chinese abacus leads to *inquiry based learning* in mathematics: asking questions, discuss different procedures (Poisard and Gueudet 2010). An important factor for teacher to do so, is the possibility to show the different procedures to all students, and is thus linked to the material environment available in the classroom (IWB, video-projector, viewer with a material abacus).

---

4   CP : *Cours Préparatoire*, CM : *Cours Moyen*.
5   In France, primary teachers can teach students aged from 2-3 years old (first year of "école maternelle") to 10-11 years old (last year of "école élémentaire").

At last, concerning the impact on classroom practices using the abacus, Ester thinks she evolved on her vision of mathematics and on her mathematical practice in the classroom. "I think before I manipulated less, that means thanks to the abacus, I went back on basics of mathematical manipulations, even with older students [year 6, 10-11 years old]. It is a change in my practice. Now as a trainer, I am more on the importance to take time, and do not make abstraction of this material appropriation, almost physical of the number. I myself modified my vision of mathematics!". Deborah locates her evolution on the "conception of number as a teacher", the attention given to make explicit different students procedures and the use of the virtual abacus in year 1 (5-6 years old) for number sense: "What was very interesting for me is the use of new technologies with young students. They are without prior assumptions, on the contrary, they easily manipulate this tool and what is really interesting as well is to see the way they go from the material abacus to the virtual one". For Rose, it was the first year of integration in her class, next year she would like to add some "collective sessions to stabilise learnings" and extend the study to calculation as well.

We next analyse the type of task: "*to set a number on the Chinese abacus*". We describe three tasks referring to different registers: abaci, fingers, and oral registers. The Chinese abacus is a dynamic resource, we show that the observation of students' gestures is important to understand their learning.

## 3. PRAXEOLOGICAL ANALYSIS OF TASKS GIVEN TO STUDENTS IN DIFFERENT REGISTERS

We analyse each task in different registers, that means we present an articulation of registers for a given task. From our observations, we identified these different registers being developed by teachers. We saw (in part 2) that in teachers' discourse, the notion of register is an underlying idea mentioned. The second task we describe ("to set the number 8") was observed in Deborah's class (Gueudet et al, 2014). We give here a praxeological analysis of tasks for students and identify knowledge related to registers.

Students' activity is a very important resource for teaching. Teacher must follow each student activity to use it for teaching. We focus on three tasks given to students: to set 3, to set 8, and to set and say 73, and we consider four registers:

- $R_A$: the Chinese abacus register, with two parts that can be distinguished (or not): $R_{VA}$ (virtual) and $R_{MA}$ (material)
- $R_W$: paper-and-pencil worksheet with drawings of abaci and beads
- $R_F$: the fingers register with two hands
- $R_O$: the oral languages register

For each task and register, we name some techniques (T) and analyse them to show the mathematical knowledge associated that is the technology ($\tau$). We do not discuss here any hierarchy of these registers, we show how they can be articulated for a specific mathematical knowledge.

We talk about *coding* numbers, for us all these registers are ways to code numbers. The term code is used here as a system of letters, numbers, symbols, etc. that represents another system, in order to codify this system. Abaci, fingers, digits, letters, drawings, cubes, etc. We consider all these resources with a common feature: the possibility to code numbers.

The Chinese abacus is an artefact that becomes an instrument when combined to students schemes of use. For this register $R_A$, the technique in an instrumented technique. The worksheets used with drawings of abacus use the same way to code numbers as the abaci register, so the technique is close to an instrumented one.

For Radford (2009), gesture is an important part of the learning of mathematical activity.

"My point here is not to diminish the cognitive role of the written. It is rather an invitation to entertain the idea that mathematical cognition is not only mediated by written symbols, but that is also mediated, in a *genuine sense*, by actions, gestures and other types of signs." (Ibid, p.112).

In our study, the fingers register is very special and involves the body to represent numbers. This register is very often used in class, and if we look carefully at it, different techniques are possible to set a number (see task 1 and 2 below) that means fingers are associated to schemes and can be considered as instruments as well. Though we could talk about instrumented technique for fingers.

The following table 1 presents the notations we use for the praxeological analysis.

| Registers R | Techniques T | Technologies τ |
|---|---|---|
| • $R_A$ (abacus): $R_{VA}$ (virtual) and $R_{MA}$ (material)<br>• $R_W$ (worksheet)<br>• $R_F$ (fingers)<br>• $R_O$ (oral languages) | For a given task and a register, the technique is the way of doing, of answering the mathematical question. The technique is noted (R, T). | The technologie is the mathematical kwnoledge associated to the technique. The technologie depends on the register as well and is noted (R, τ). |

**TABLE 1:** Notations used for the praxeological analysis of tasks in different registers

*3.1 Task 1: "to set the number 3"*
For this task only beads in the units are activated.

- $R_A$ abaci register

To active a bead on the material abacus means to move the beads, up for 1-unit counters and down for the 5-unit counters. The activation of a bead on the virtual abacus is made by clicking on a bead. To active three, it is possible to click on the third one. When techniques are similar on both the material and the virtual abacus, we do not need to separate the two registers $R_{MA}$ and $R_{VA}$ (see $T_1$ and $T_2$).

| Techniques T | Technologies τ |
|---|---|
| ($R_A$, $T_1$)<br>To activate one 1-unit counter in the units, three times one after the other, that means in three gestures. Like saying: 1, 2, 3.<br>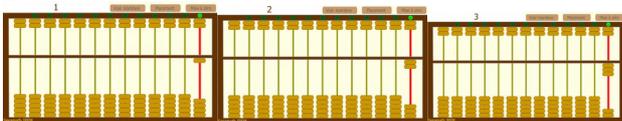 | ($R_A$, $τ_1$)<br>This is a *counting reasoning:*<br>1, 2, 3.<br>Students know the numbers sequence maybe just like a rhyme. |
| ($R_A$, $T_2$)<br>To activate the third bead in the units, in one gesture.<br>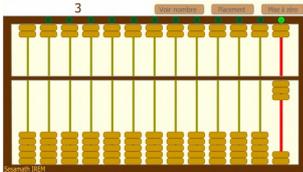 | ($R_A$, $τ_2$)<br>This is a reasoning referring to the ordinality of number: *ordinality reasoning:*<br>the third.<br>Students can identify the third bead without counting the three beads. |
| ($Rv_A$, $T_3$)<br>On the virtual abacus, written number can be | ($Rv_A$, $τ_3$)<br>This is a *trial/error approach*. |

shown and students can try different possibilities to display the number 3 on the software (upper left).

- $R_W$ paper-and-pencil worksheet register

| Techniques T | Technologies τ |
|---|---|
| ($R_W$, $T_1$) To draw all beads, activated or not 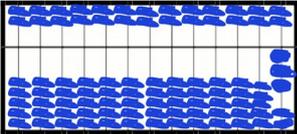 | ($R_W$, $τ_1$) This student needs to represent all beads on this abacus. She/he sees this abacus as an artefact and may be far from a mathematical meaning. |
| ($R_W$, $T_2$) To draw the three activated beads 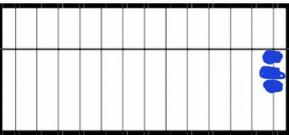 | ($R_W$, $τ_2$) This is a first abstraction: this student represents only the three activated beads that give the meaning of three. |
| ($R_W$, $T_3$) To symbolise the three activated beads with a line. 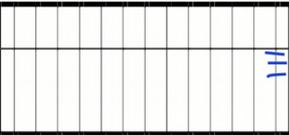 | ($R_W$, $τ_3$) As this drawing has some abstraction, (and is not a drawing of all beads), most students using this technique have a good understanding of the abacus and the way it works. They have reached a mathematical work. |

- $R_F$ fingers register

| Techniques T | Technologies τ |
|---|---|
| ($R_F$, $T_1$) To set three fingers on one hand. This technique T1 has several under-techniques: T1a, T1b, T1c, depending on the fingers: from the little one of from the thumb (like we do in France). 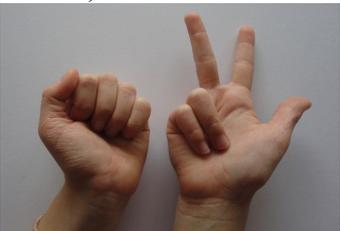 | ($R_F$, $τ_1$) Usually, only one hand is used in France to show 3. With two hands, we show that: 3=0+3 |
| ($R_F$, $T_2$) To set two fingers on one hand and one finger on the other one. This technique T2 has as well several under-techniques: T2a, T1b, T2c, depending on the fingers. | ($R_F$, $τ_2$) With two hands, we can show another decomposition of 3: 3=1+2 |

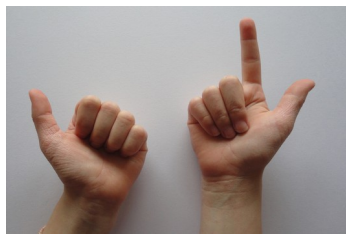

*3.2 Task 2: "to set the number 8"*
For this task only beads in the units are activated.

- $R_A$ abaci register: material ($R_{MA}$) and virtual ($R_{VA}$)

For the three first techniques, they are possible on both the material and the virtual Chinese abacus. For the fourth one, the technique is available on the material but not on the virtual one.

| Techniques T | Technologies τ |
|---|---|
| ($R_A$, $T_1$)<br>To activate one 5-unit counter.<br>To activate three 1-unit counters.<br>That means two gestures.<br>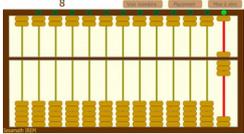 | ($R_A$, $τ_1$)<br>This is a reasoning that reveals that this student knows well to decompose numbers additively, we call it a *calculating reasoning*.<br>8=5+3 (two gestures) |
| ($R_A$, $T_2$)<br>To activate one 5-unit counter first. After to activate three 1-unit counters in three gestures. | ($R_A$, $τ_2$)<br>This student knows that the beads in the upper part count for five. She/he is able to make the difference between *quantity* (of beads) and *value* (of beads here five per rod). After five, this student over-count, one by one. She/he is more confident up to five and after needs to over-count.<br>8=5+1+1+1 |
| ($R_A$, $T_3$)<br>To activate five 1-unit counters, in five gesture.<br>To exchange five 1-unit counters against one 5-unit counter (in the units).<br>To activate three 1-unit counters in three gestures. | ($R_A$, $τ_3$)<br>This is a *counting reasoning*.<br>8=(1+1+1+1+1)+1+1+1 and 1+1+1+1+1=5 |
| ($R_A$, $T_4$)<br>To activate five 1-unit counters, in five gesture.<br>To exchange five 1-unit counters against one 5-unit counter (in the units).<br>To activate the third 1-unit counter in one gesture. | ($R_A$, $τ_4$)<br>This student has a *counting reasoning* to set five and an *ordinality reasoning* for three. She/he is more confident with a small number like three and needs to make the exchange to set five.<br>8=(1+1+1+1+1)+3 |
| ($R_{MA}$, $T_5$): material abacus<br>To activate one 5-unit counter and three 1-unit counters (in the units) in one gesture.<br>This technique is available on the material abacus but not on the virtual one we use[6].<br>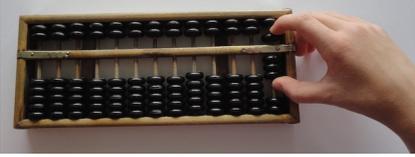 | ($R_{MA}$, $T_5$)<br>This is a reasoning that revels that this student knows well to decompose numbers, we call it a *calculating reasoning*.<br>8=5+3 (one gesture) |
| ($Rv_A$, $T_6$)<br>On the virtual abacus, written number can be shown and students can try different possibilities to display the number 8 on the software (upper left). | ($Rv_A$, $\tau_6$)<br>This is a *trial/error approach*. |

- $R_F$ fingers register

These techniques have several under-techniques depending on the fingers that are shown. Some are easier due to our physical particularity but as well due to the usual way to set numbers in a given culture (see task 1). This register is usual in classrooms and often used for counting with particularities depending on the culture. In France, we always start with the thumb to count on fingers for example. If we see someone starting with the little finger or the index finger, we know she/he is not French! For example, different studies looked at African languages and gestures about counting and had shown that counting on fingers is not universal but is a cultural and social result (Zaslavsky 1973, Gerdes 2009). The chambaa and the makonde ways of counting on fingers is very interesting to study (Gerdes, 2009, p.133-134). The chambaa people use one hand to set four but grouping twice two (4=2+2, one hand), six is set on two hands (6=3+3), for seven fingers show 7=(2+2)+3 and eight is 8=(2+2)+(2+2). Another example: the makonde people use two hands to count up to four and the folded fingers mark the numbers: one is the folded little finger and one more is folded up to four.

We show here the task to set a number on fingers has different techniques linked to the different decomposition of numbers and is a resource for teaching. In this analysis we show only two techniques. We think we can use fingers to show some other decompositions of 8 but it is delicate to describe it here as it has a dynamic way of grouping fingers. The decomposition 8=6+2 on hands can be: first 8=5+3 and after the thumb of the second hand in the fist of the other hand. Similar gestures are possible to show that 8=8+0 and 8=1+7.

| Techniques T | Technologies τ |
|---|---|
| ($R_F$, $T_1$)<br>To set five fingers on one hand, and three fingers on the other one.<br>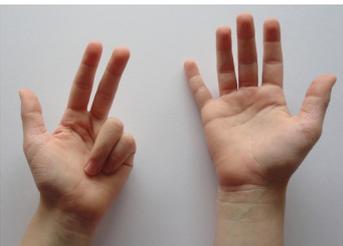 | ($R_F$, $\tau_1$)<br>This technique is very often used, sometimes it is the only one used by primary teachers in this register (in France).<br>8=3+5 |

---
6  Some virtual abacus, in particular applications for tablets can set eight in one gesture.

| (R_F, T_2) To set four fingers on one hand, and four fingers on the other one. 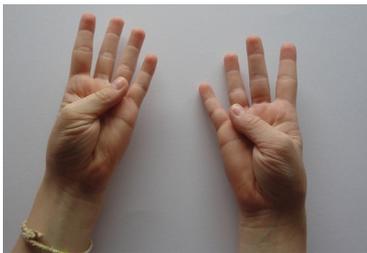 | (R_F, τ_2) This technique work on *addition of doubles* very important to learn calculation: 8=4+4 |
|---|---|

*3.3 Task 3: "to set and say the number 73"*

- $R_O$ oral languages register

A language reflects the way we see the world and mathematics is a modelling of the world. Each language has a special way to say numbers, space, measuring, etc. (see for example Barton, 2006). The mathematical modelling -that is to say the technology- associated to an oral number system is different for every language. We show that for the example of the way numbers are said in four languages[7]: English, French, Maori and Breton. (In the Breton case, see Poisard et al, 2014 and Kervran et al 2015 for more details). We present here languages as a resource for the teaching and learning of mathematics, and more precisely as both a register and a technique. Indeed, we argue that looking at other languages helps understanding our own language (or languages); and for the mathematical vocabulary helps understanding mathematics. We take here four very different languages for oral number system (for some other examples, see Ascher 2002).

| Techniques T | Technologies τ |
|---|---|
| (R_O, T_1) in English Seventy-three | (R_O, τ_1) The English 73 is not far from a regular form that would be *seven-ten (and) three*: 73=7×10+3. |
| (R_O, T_2) in French *soixante-treize* literally *sixty-thirteen* | (R_O, τ_2) In French, there is a reference to a 60 grouping that we can find in time (60 seconds in a minute): 73=60+13 |
| (R_O, T_3) in Maori *whitu tekau ma toru* literally *seven-tens and three* | (R_O, τ_3) Modern Maori has a regular way to say numbers that is to say the way numbers are said is the same as the mathematical meaning (polynomial development): 73=7×10+3 |
| (R_O, T_4) in Breton *trizek ha tri-ugent,* literally *three-ten and three-twenty.* (*thirteen* is said *three-ten* in Breton). | (R_O, τ_4) As a Celtic language, Breton use vigesimal (base 20) grouping to say numbers: 73=3+10+3×20. |

- $R_A$ abaci register

---

7  Acknowledgments to Tony Trinick (Maori) and Erwan Le Pipec (Breton).

| Techniques T | Technologies τ |
|---|---|
| ($R_A$, $T_1$)<br>In the units: to activate three 1-unit counters.<br>In the tens: to activate one 5-unit counter and two 1-unit counters.<br>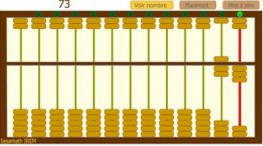 | ($R_A$, $τ_1$)<br>This is the *economical inscription* of 73:<br>This is very similar to what is said in Maori or in Chinese.<br>$73 = 7 \times 10 + 3$ |
| ($R_A$, $T_2$)<br>In the units: to activate three 1-unit counters and two 5-unit counters.<br>In the tens: to activate one 5-unit counter and one 1-unit counter.<br>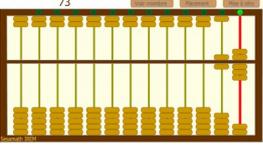 | ($R_A$, $τ_2$)<br>This is similar to the French way of saying 73: *soixante-treize* (literally *sixty-thirteen*).<br>$73 = 60 + 13$ |
| ($Rv_A$, $T_3$)<br>On the virtual abacus, written number can be shown and students can try different possibilities to display the number 3 on the software (upper left). | ($Rv_A$, $τ_3$)<br>This is a *trial/error approach*. |

The analysis of these three tasks referring to four different registers show that the mathematical knowledge depends on both the register and the technique used by students. Depending on students' choices of technique, we analyse her/his knowledge. The mathematical knowledge describes here is about the difference between quantity and value, the counting, ordinality, and calculating reasonings, the decomposition of numbers. The notion of commutativity of the addition is also raised. For one register, different techniques are available, and we can see that sometimes one technique is possible in different registers (Table 2).

| *To set 3* | $R_A$ Abaci | $R_W$ Worksheet | $R_F$ Fingers |
|---|---|---|---|
| | $(R_A, T_1, \tau_1)$: counting | $(R_W, T_1, \tau_1)$: all drawings | $(R_F, T_1, \tau_1)$: 3=0+3 |
| | $(R_A, T_2, \tau_2)$: ordinality | $(R_W, T_2, \tau_2)$: drawings | $(R_F, T_1, \tau_1)$: 3=1+2 |
| | $(Rv_A, T_3, \tau_3)$: trial/error | $(R_W, T_3, \tau_3)$: symbols | |

| *To set 8* | $R_A$ Abaci | $R_F$ Fingers |
|---|---|---|
| | $(R_A, T_1, \tau_1)$: 8=5+3 (two gestures), calculating | $(R_F, T_1, \tau_1)$: 8=3+5 |
| | $(R_A, T_2, \tau_2)$: 8=5+1+1+1 and quantity/value | $(R_F, T_2, \tau_2)$: 8=4+4 |
| | $(R_A, T_3, \tau_3)$: 8=(1+1+1+1+1)+1+1+1 and 1+1+1+1+1=5, counting | |
| | $(R_A, T_4, \tau_4)$: 8=(1+1+1+1+1)+3, counting and ordinality | |
| | $(R_{MA}, T_5, \tau_5)$: 8=5+3 (two gestures), calculating | |
| | $(Rv_A, T_6, \tau_6)$: trial/error | |

| *To set and say 73* | $R_A$ Abaci | $R_O$ Oral |
|---|---|---|
| | $(R_A, T_1, \tau_1)$: economical inscription, 73=7×10+3 | $(R_A, T_1, \tau_1)$: English, 73=7×10+3 |
| | $(R_A, T_2, \tau_2)$: 73=60+13 | $(R_A, T_2, \tau_2)$: French, 73=60+13 |
| | $(Rv_A, T_3, \tau_3)$: trial/error | $(R_A, T_3, \tau_3)$: Maori, 73=7×10+3 |
| | | $(R_A, T4, \tau4)$: Breton, 73=3+10+3×20 |

**TABLE 2:** Articulation of registers for the three tasks

The articulation of these registers are very important for teachers. A teacher can give a task orally or written, and this is different to say a number (see the example of 73), to write it in digits, or in letters as well. From this task, students have to set it on the abacus, the material one, the virtual one, or a worksheet but it is also possible on fingers. The virtual abacus has a very specific feature as the number is set on the abacus and can be shown in numbers as well. The abacus is introduced as a resource for teaching but can not replace some other resources. We can see it clearly with the task *to set 8*: the Chinese abacus allows to work on the decomposition 5+3 which is very important to learn mental calculation, but does not work on other decompositions of 8 like 4+4. Fingers allow to work on 4+4 for example. Looking carefully at students' activity and gestures to set 8 can help teacher to understand students' level about: value and quantity, counting up to 5 and 8, ordinality, and calculating. Moreover, oral number system is a resource for teaching, a language is not transparent to teach mathematics. The analysis of the task *to set and say 73* shows that the abacus can help understand the French way of saying 73 (sixty-thirteen), whereas the English and Maori languages are close to the economical inscription. In the Breton case, we think that to identify clearly for students the technology: 73=3+10+3×20 is a resource for teaching Breton oral number system.

## CONCLUSION
This work is about coding numbers and the associated knowledge about this *code*. For example, for a small number like three, we can say it, write it in digits, write it in letters, show it on one hand, show it on two hands, set it on an abacus, say it in another language, etc. All these registers refer to a specific knowledge, more precisely for one register different techniques are possible and a technique can be used in different registers. This is this articulation between registers and techniques that makes it a resource for teaching and learning mathematics. Students need to recognise a number in different registers and techniques to understand the meaning of what is a number. We see it described by Ester, Deborah and Rose in their interviews and we analyse it precisely for three tasks. We show how the articulation between registers and techniques is possible. For Duval, if a student is able to achieve a task in at least two registers, it means this student understands the mathematical knowledge concerns. We think that to identify students' understanding the articulation between a register and the notion of technique/technology is more precise. For a given task, we argue that:

- for one register, at least two techniques are needed to make sure of students' understanding
- for one technique, at least two registers are needed to make sure of students' understanding

Identifying students' activity is central in teachers' documentation work. This is a central resource for teaching. And this is as well a central question for teacher training: the description and the analysis of different students' techniques (called sometimes procedures) and registers should have a large place in trainings.

Moreover, this study makes an important place to new technologies: the virtual abacus is used on computers by both students and teachers, the video-projector and the IWB is a central resource for teachers in the class, etc. With the material abacus, some teachers use a viewer, so the material abacus is entering in the digital area as well! But in some schools, with less technologies a paper poster of the abacus with magnets makes it possible to be used in class. We show here that new technologies bring new appropriations by teachers and this is also a question of articulation: between material and virtual resources.

**REFERENCES**

Adler, J. (2000). Conceptualising resources as a theme for teacher education. *Journal of Mathematics Teacher Education. 3. 205-224.*
Artigue, M. (2009). Rapports et articulations entre cadres théoriques: le cas de la théorie anthropologique du didactique. *Recherches en didactique des mathématiques. 29(3)*. 305-334.
Asher, M. (2002). *Mathematics Elsewhere: An Exploration of Ideas Across Cultures.* Princeton: Princeton University Press.
Baccaglini-Frank, A., & Maracci, M. (2015). Multi-Touch Technology and Preschoolers' Development of Number-Sense. *Digital Experiences in Mathematics Education*, *1*(1), 7-27.
Bartolini, M. G., & Martignone, F. (2014). Manipulatives in Mathematics Education. In S. Lerman (Éd.), *Encyclopedia of Mathematics Education* (p. 365-372). Dordrecht: Springer Netherlands.
Barton B. (2008) *The language of mathematics. Telling mathematical tales.* New York : Springer.
Bueno-Ravel, L. & Harel, C. (2016, to appear). Le calcul mental à l'école: apports du boulier chinois. *Mathematice 51.*
Chevallard, Y. (2007) Readjusting didactics to a changing epistemology. *Eur Educ Res J 6(2):*131–134
Chevallard, Y. (2006) Steps towards a new epistemology in mathematics education. In: Bosch M (ed) *Proceedings of CERME 4.* Fundemi IQS, Barcelone, pp 21–30
Chevallard, Y. (1999). L'analyses des pratiques enseignantes en théorie anthropologique du didactique. *Recherches en Didactique des Mathématiques*, 19(2), 221-266.
Chevallard Y (1992b) Fundamental concepts in didactics: perspectives provided by an anthropological approach. In: Douady R, Mercier A (eds) *Research in didactique of mathematics, selected papers.* La pensée sauvage, Grenoble. p.131-167
Chevallard, Y. (1992a). Concepts fondamentaux de la didactique: Perspectives apportées par une approche anthropologique. *Recherches en Didactique des Mathématiques*, 12(1), 73-112.
Duval, R. (2006). A Cognitive Analysis of Problems of Comprehension in a Learning of Mathematics. *Educational Studies in Mathematics*, *61*(1), 103-131.

Duval, R. (1998). Signe et objet (I): trois grandes étapes dans la problématique des rapports entre représentation et objet. *Annales de Didactique et de Sciences Cognitives* 6, 139-163.

Duval, R. (1996). Quel cognitif retenir en didactique des mathématiques? *Recherches en Didactique des Mathématiques* 16(3), 349-382.

Gerdes, P. (2009). *L'ethnomathématique en Afrique.* Maputo: CEMEC.

Gueudet, G., Bueno-Ravel, L., & Poisard, C. (2014). Teaching Mathematics with Technology at the Kindergarten Level: Resources and Orchestrations. In A. Clark-Wilson, O. Robutti, & N. Sinclair (Éd.), *The Mathematics Teacher in the Digital Era: An International Perspective on Technology Focused Professional Development* (p. 213-240). Dordrecht: Springer Netherlands.

Gueudet, G., & Trouche, L. (2009). Towards new documentation systems for mathematics teachers? *Educational Studies in Mathematics, 71*(3), 199–218.

Guin, D., Ruthven, K., & Trouche, L. (2005). *The Didactical Challenge of Symbolic Calculators: Turning a Computational Device into a Mathematical Instrument.* Boston, MA: Springer US.

Kervran, M., Poisard, C., Le Pipec, E, Sichler, M., & Jeudy-Karakoç, N. (2015). Langues minoritaires locales et conceptualisation à l'école : l 'exemple de l'enseignement des mathématiques en breton. In Kervran M., Blanchet, P. (dir). *Langues minoritaires locales et éducation a la diversité des dispositifs didactiques a l'épreuve.* L'harmattan, coll. Espaces Discursif, 65-82.

Lagrange, J.-B. (2005). Using Symbolic Calculators to Study Mathematics. In D. Guin, K. Ruthven, & L. Trouche (Éd.), *The Didactical Challenge of Symbolic Calculators: Turning a Computational Device into a Mathematical Instrument* (p. 113-135). Boston, MA: Springer US.

Maschietto, M., & Soury-Lavergne, S. (2013). Designing a duo of material and digital artifacts: the pascaline and Cabri Elem e-books in primary school mathematics. *ZDM*, *45*(7), 959-971.

Maschietto, M., & Trouche, L. (2010). Mathematics learning and tools from theoretical, historical and practical points of view: the productive notion of mathematics laboratories. *ZDM*, *42*(1), 33-47.

Block, D., Nikolantonakis, K., & Vivier, L. (2012). Registres et praxis numérique en fin de première année de primaire dans trois pays. *Annales de didactique et de sciences cognitives. 17*. 59-86.

Poisard, C. (2005b). Ateliers de fabrication et d'étude d'objets mathématiques, le cas des instruments à calculer. Thèse de Doctorat de l'Université de Provence, Aix-Marseille I.

Poisard, C., Gueudet, G., & Bueno-Ravel, L. (2011). Comprendre l'intégration de ressources technologiques en mathématiques par des professeurs des écoles. *Recherches en Didactique des Mathématiques*, *31(2)*, 151-189.

Poisard, C., & Gueudet, G. (2010). Démarches d'investigation : exemples avec le boulier virtuel, la calculatrice et le TBI. *Journées mathématiques de l'INRP*, Lyon.

Poisard, C., Kervran, M., Le Pipec, E., Alliot, S., Gueudet, G., Hili, H., Jeudy-Karadoc, N., & Larvol, G. (2014). Enseignement et apprentissage des mathématiques à l'école primaire dans un contexte bilingue breton-français. *Revue Spirale 54. 129-150.*

Poisard, C. (2005a). Les objets mathématiques matériels, l'exemple du boulier chinois, *Petit x*, 68, 39-67.

Poisard, C., Riou-Azou, G., D'hondt, D., & Moumin, E. (2016, to appear). Le boulier chinois : une ressource pour la classe et pour la formation des professeurs. *Mathematice 51.*

Radford, L. (2009). Why do gestures matter? Sensuous cognition and the palpability of mathematical meanings. *Educational Studies in Mathematics*, *70*(2), 111-126.

Rau, P.-L. P., Xie, A., Li, Z., & Chen, C. (2016). The Cognitive Process of Chinese Abacus Arithmetic. *International Journal of Science and Mathematics Education*, *14*(8), 1499-1516.

Riou-Azou, G. (2014). La construction du nombre en grande section de maternelle avec un boulier chinois virtuel. *Mathematice 40.*

Sun, X., Berinderjeet, K. & Novotna, J. (2015). *Primary Mathematics Study on Whole numbers.* Conference proceedings of ICMI Study 23.

Trouche, L. (2005). An Instrumental Approach to Mathematics Learning in Symbolic Calculator Environments. In D. Guin, K. Ruthven, & L. Trouche (Éd.), *The Didactical Challenge of Symbolic Calculators: Turning a Computational Device into a Mathematical Instrument* (p. 137-162). Boston, MA: Springer US.

Verillon, P., & Rabardel, P. (1995). Cognition and Artefacts: A contribution to the study of thought in relation to instrumented activity. *European Journal of Psychology of Education, 10*, 77–102.

Zaslavsky C. (1973, 1999). *Africa counts. Number and pattern in african cultures. Third edition.* Chicago: Lauwrence Hill Books.